\documentclass[english,conference]{ieeeconf}
\usepackage[T1]{fontenc}
\usepackage[latin9]{inputenc}
\usepackage{mathtools}
\usepackage{amsmath}
\usepackage{amsthm}
\usepackage{amssymb}
\usepackage{graphicx}

\makeatletter
\theoremstyle{plain}
\newtheorem{thm}{\protect\theoremname}
\theoremstyle{definition}
\newtheorem{example}[thm]{\protect\examplename}
\theoremstyle{remark}
\newtheorem{rem}[thm]{\protect\remarkname}

\usepackage[caption=false,font=footnotesize]{subfig}
\usepackage{url}
\usepackage[nocompress]{cite}

\usepackage{amsmath}
\usepackage{amssymb}
\usepackage{amsfonts}
\usepackage{mathtools}
\usepackage{autonum}

\newcommand{\restore@Environment}[1]{%
	\AtBeginDocument{%
		\csletcs{#1*}{#1}%
		\csletcs{end#1*}{end#1}%
	}%
}
\forcsvlist\restore@Environment{alignat,equation,gather,multline,flalign,align}

\usepackage{arydshln}

\newcommand{\norm}[1]{\left\Vert#1\right\Vert}

\newcommand{\reals}{\mathbb{R}}

\newcommand{\optmin}{\mathrm{min.}}

\newcommand{\vOne}{\mathbf{1}}

\newcommand{\xopt}{x^{\star}}

\newcommand{\vopt}{v^{\star}}
\newcommand{\zopt}{z^{\star}}
\newcommand{\wopt}{w^{\star}}
\newcommand{\opt}{^{\star}}

\newcommand{\deq}{\coloneqq}

\newcommand{\Gcon}{G_\mathrm{con}}
\newcommand{\ave}{\mathsf{ave}}

\newcommand{\Jo}{J_\perp}

\newcommand{\Gadmm}{G_\mathrm{ADMM}}

\allowdisplaybreaks[3]
\IEEEoverridecommandlockouts
\newcommand{\allthanks}{\thanks{The author is with the Department of Electrical and Computer Engineering, University of Illinois, Chicago, IL 60607. {\tt\scriptsize {hanshuo@uic.edu}}.}}

\@ifundefined{showcaptionsetup}{}{%
 \PassOptionsToPackage{caption=false}{subfig}}
\usepackage{subfig}
\makeatother

\usepackage{babel}
\providecommand{\examplename}{Example}
\providecommand{\remarkname}{Remark}
\providecommand{\theoremname}{Theorem}

\begin{document}

\title{\textbf{\Large{}Systematic Design of Decentralized Algorithms for Consensus Optimization}}

\author{Shuo Han\allthanks}
\maketitle
\begin{abstract}
We propose a separation principle that enables a systematic way of designing decentralized algorithms used in consensus optimization. Specifically, we show that a decentralized optimization algorithm can be constructed by combining a non-decentralized base optimization algorithm and decentralized consensus tracking. The separation principle provides modularity in both the design and analysis of algorithms under an automated convergence analysis framework using integral quadratic constraints (IQCs). We show that consensus tracking can be incorporated into the IQC-based analysis. The workflow is illustrated through the design and analysis of a decentralized algorithm based on the alternating direction method of multipliers. 

\end{abstract}

\section{Introduction}

In this paper, we study algorithms for solving the \emph{consensus optimization problem}, which has the form
\begin{equation}
\underset{x_{0}\in\mathbb{R}^{d}}{\optmin}\quad f_{0}(x_{0})\coloneqq\frac{1}{n}\sum_{i=1}^{n}f_{i}(x_{0}).\label{eq:prob}
\end{equation}
We assume $f_{i}\colon\mathbb{R}^{d}\to\mathbb{R}$ is convex for $i\in\{1,2,\dots,n\}$ and the set of minimizers is nonempty. The name consensus optimization is due to the fact that the problem can be made equivalent to another optimization problem with a separable objective function $\sum_{i=1}^{n}f_{i}(x_{i})$ by introducing local optimization variables $x_{1},x_{2},\dots,x_{n}\in\mathbb{R}^{d}$ and a consensus constraint $x_{1}=x_{2}=\dots=x_{n}$. 

We are interested in algorithms that solve the consensus optimization problem in a \emph{decentralized} manner. We shall make a distinction between distributed and decentralized algorithms, which are often used interchangeably in the literature; the former permits the presence of a \emph{master node} that collects computational results from multiple worker nodes, whereas the latter does not require a master node. Most existing decentralized algorithms used in consensus optimization belong to one of the following two classes. The first one is based on the gradient descent method or its variants (e.g., Nesterov's method). This includes, among others, the distributed gradient descent method\ \cite{nedic_distributed_2009} (which is, in fact, decentralized despite its name), DIGing\ \cite{nedic_achieving_2017,qu_harnessing_2018}, and EXTRA\ \cite{shi_extra:_2015}. See also\ \cite{qu_accelerated_2017} for an algorithm based on Nesterov's method and\ \cite{xin_linear_2018} for handling directed communication graphs. The second one is based on operator splitting methods, of which the most widely used is the Douglas\textendash Rachford method\ \cite{parikh_proximal_2014} or its application to the dual problem, the alternating direction method of multipliers (ADMM)\ \cite{boyd_distributed_2011}. Although the original ADMM algorithm, when directly applied to the consensus optimization problem, requires a master node and therefore is not decentralized, it has been shown that ADMM can be made decentralized through a reformulation of the consensus constraint\ \cite{wei_o1/k_2013,shi_linear_2014}. 

Despite a vast body of literature on decentralized optimization algorithms in recent years, there has been little work on systematic understanding and designing of decentralized algorithms. (See\ \cite{sundararajan_canonical_2018} for some recent effort on unifying decentralized algorithms that are based on gradient descent.) As a result, whenever the base optimization algorithm changes (e.g., from regular gradient descent to accelerated gradient descent) or the conditions on the communication graph changes (e.g., from undirected to directed), a convergence analysis of the new algorithm needs to be started almost from scratch. This paper seeks a framework that enables a systematic design of decentralized optimization algorithms in the hope of speeding up the development of new algorithms. 

We believe such a framework can be made possible through the automated convergence analysis of optimization algorithms proposed recently by Lessard et al.\ \cite{lessard_analysis_2016} Unlike traditional, proof-based analysis that needs to be carried out manually, their automated convergence analysis uses computational tools to establish a numerical certificate of convergence for optimization algorithms. The key is to view an optimization algorithm as a feedback interconnection of a linear dynamical system and a nonlinear memoryless but uncertain system, which can be characterized by \emph{integral quadratic constraints} (IQCs)\ \cite{megretski_system_1997}. As a result, convergence of optimization algorithms can be established by certifying stability of the feedback interconnection.

\paragraph*{Contribution}

The main contribution of this paper is a \emph{separation principle} for the design of decentralized algorithms used in consensus optimization. Specifically, one can start with a non-decentralized optimization algorithm and replace the (static) averaging operation therein with decentralized average consensus tracking. Such an approach not only enables a systematic way for designing decentralized optimization algorithms but is also amenable to automated convergence analysis based on IQC. We believe that the result will help unify existing decentralized algorithms and eventually facilitate the development of new algorithms.

When applied to known settings (i.e., same base algorithm, same conditions on the communication graph), the result of this paper is not guaranteed to yield a better (e.g., with faster convergence, or more robust) decentralized algorithm than existing ones; our main focus is a more principled design procedure rather than optimality. The use of consensus tracking in decentralized optimization is not new and can be found in the DIGing algorithm\ \cite{nedic_achieving_2017,qu_harnessing_2018}; the role of this paper is to highlight how consensus tracking can be separated from the base algorithm.

\section{Main Results~\label{sec:results}}

\subsection{Notation}

Denote by $\vOne$ the column vector of all ones, $\norm{\cdot}$ the $\ell_{2}$-norm of a vector, $I_{n}$ the $n\times n$ identity matrix (size omitted when clear from the context), $\otimes$ the Kronecker product, and $\sigma_{\max}(\cdot)$ the maximum singular value of a matrix. We also define $J\deq\vOne\vOne^{T}/n$ and $\Jo\deq I-J$. We use exclusively $W$ to denote a symmetric irreducible doubly stochastic matrix (called \emph{gossip matrix} in the setting of consensus): $W\vOne=\vOne$ and $\vOne^{T}W=\vOne^{T}$. For a symmetric matrix $P$, we write $P\succeq0$ if $P$ is positive semidefinite. For a differentiable function $f$, we denote by $\nabla f$ the gradient of $f$. 

In the context of decentralized optimization, each node is required to keep its own local variables. We reserve the subscript for indexing the nodes and the superscript for indexing a given sequence. For example, $x_{i}\in\mathbb{R}^{d}$ represents a local variable that belongs to node $i$, whereas a sequence of vectors is denoted by $\{x^{k}\}_{k\geq0}\coloneqq\{x^{0},x^{1},\dots\}$. For any convergent sequence $\{x^{k}\}$, we use $\xopt$ to denote its limit or, alternatively, steady-state value. We use the notation 
\[
x\deq\left[\begin{array}{cccc}
x_{1} & x_{2} & \cdots & x_{n}\end{array}\right]^{T}\in\reals^{n\times d}
\]
to denote the matrix whose rows are formed by local variables $x_{1},x_{2},\dots,x_{n}\in\reals^{d}$, and we define $\ave(x)\coloneqq\frac{1}{n}\sum_{i=1}^{n}x_{i}=\frac{1}{n}x^{T}\vOne\in\mathbb{R}^{d}$. Similarly, we use the notation 
\[
\nabla f(x)\deq\left[\begin{array}{ccc}
\nabla f_{1}(x_{1}) & \cdots & \nabla f_{n}(x_{n})\end{array}\right]^{T}\in\reals^{n\times d}
\]
to denote the list of local gradients. 

\subsection{Problem description}

We investigate algorithms that solve the consensus optimization problem and can be expressed in one of the following two forms.
\begin{enumerate}
\item Centralized algorithms:\begin{subequations}\label{eq:calg}
\begin{align}
\xi_{0}^{k+1} & =A_{0}\xi_{0}^{k}+B_{0}\,\ave(u^{k})\label{eq:calg-a}\\
v_{0}^{k} & =C_{0}\xi_{0}^{k}+D_{0}\,\ave(u^{k})\label{eq:calg-b}\\
u_{i}^{k} & =\phi_{i}(v_{0}^{k}),\quad i=1,2,\dots,n.\label{eq:calg-c}
\end{align}
\end{subequations}
\item Distributed algorithms: for $i\in\{1,2,\dots,n\}$, \begin{subequations}\label{eq:dalg}
\begin{align}
\xi_{i}^{k+1} & =A_{i}\xi_{i}^{k}+B_{\text{loc},i}u_{i}^{k}+B_{i}\,\ave(u^{k})\label{eq:dalg-a}\\
v_{i}^{k} & =C_{i}\xi_{i}^{k}+D_{\text{loc},i}u_{i}^{k}+D_{i}\,\ave(u^{k})\label{eq:dalg-b}\\
u_{i}^{k} & =\phi_{i}(v_{i}^{k}).\label{eq:dalg-c}
\end{align}
\end{subequations}
\end{enumerate}
In both cases, $\phi_{i}$ is a continuous but possibly nonlinear function; $A_{i}$, $B_{i}$, $C_{i}$, $D_{i}$, $B_{\text{loc},i}$, and $D_{\text{loc},i}$ are all constant matrices of appropriate dimensions. (The subscript ``$\text{loc}$'' stands for \emph{local}.) We give for each case one example algorithm that solves the consensus optimization problem. To simplify notation, we assume in the remaining part of this section that each $f_{i}$ is smooth, but we expect the result to generalize to nonsmooth objective functions by making use of subdifferentials.
\begin{example}
[Gradient descent] \label{exa:grad-ori}When applied to consensus optimization, the gradient descent algorithm becomes
\[
x_{0}^{k+1}=x_{0}^{k}-\eta\nabla f_{0}(x_{0}^{k})=x_{0}^{k}-\frac{\eta}{n}\sum_{i=1}^{n}\nabla f_{i}(x_{0}^{k}),
\]
where $\eta>0$ is a constant. Define $\xi_{0}^{k}\deq x_{0}^{k}$, $v_{0}^{k}\deq\xi_{0}^{k}$, and $u_{i}^{k}\deq\nabla f_{i}(v_{0}^{k})$. The gradient descent algorithm can be written as
\[
\xi_{0}^{k+1}=\xi_{0}^{k}-\eta\,\ave(u^{k}),\qquad v_{0}^{k}=\xi_{0}^{k}.
\]
\end{example}
\begin{example}
[ADMM] \label{exa:admm-ori}When applied to consensus optimization, the ADMM algorithm becomes\ \cite[p. 50]{boyd_distributed_2011}
\begin{align*}
x_{i}^{k+1} & =\ave(x^{k})-(y_{i}^{k}+w_{i}^{k})/\rho\\
y_{i}^{k+1} & =\ave(w^{k})-w_{i}^{k}\\
w_{i}^{k} & =\nabla f_{i}(x_{i}^{k+1}),
\end{align*}
where $\rho>0$ is a constant. Define $\xi_{i}^{k}\deq(x_{i}^{k},y_{i}^{k})$, $v_{i}^{k}=(v_{i}^{1,k},v_{i}^{2,k})\deq(x_{i}^{k},x_{i}^{k+1})$, and $u_{i}^{k}\deq(v_{i}^{1,k},\nabla f_{i}(v_{i}^{2,k}))$. The ADMM algorithm can be written as
\begin{align*}
\xi_{i}^{k+1} & =\left[\begin{array}{cc}
0 & -\frac{I}{\rho}\\
0 & 0
\end{array}\right]\xi_{i}^{k}+\left[\begin{array}{c}
-\frac{I}{\rho}\\
-I
\end{array}\right]u_{i}^{2,k}+\left[\begin{array}{c}
\ave(u^{1,k})\\
\ave(u^{2,k})
\end{array}\right]\\
v_{i}^{k} & =\left[\begin{array}{cc}
I & 0\\
0 & -\frac{I}{\rho}
\end{array}\right]\xi_{i}^{k}+\left[\begin{array}{c}
0\\
-\frac{I}{\rho}
\end{array}\right]u_{i}^{2,k}+\left[\begin{array}{c}
0\\
\ave(u^{1,k})
\end{array}\right].
\end{align*}
\end{example}
The algorithms given in\ (\ref{eq:calg}) and\ (\ref{eq:dalg}) are not fully decentralized because evaluating $\ave$ requires a master node to collect information from all the nodes. Moreover, in the first case of centralized algorithms, the computation in\ (\ref{eq:calg-a}) and (\ref{eq:calg-b}) needs to be completed by the master node as well and therefore is not decentralized either. Our goal in this paper is to develop a systematic procedure for converting an existing algorithm of the form\ (\ref{eq:calg}) or\ (\ref{eq:dalg}) into a decentralized algorithm.

\subsection{Main results}

The key component in our procedure of decentralization is consensus tracking. We say that a dynamical system $\Gcon$ achieves \emph{average consensus tracking} (or simply \emph{consensus tracking}) if for any sequence $v=\{v^{k}\in\reals^{n\times d}\}$ converging to $\vopt$, the output $w=\Gcon v$ converges to $J\vopt$. An example of a system that achieves consensus tracking is given by
\begin{equation}
w^{k+1}=Ww^{k}+(v^{k+1}-v^{k}),\qquad w^{0}=v^{0}.\label{eq:gcon-1}
\end{equation}
Systems that achieve consensus tracking are not unique. For example, the system in\ (\ref{eq:gcon-1}) can be modified slightly as
\begin{equation}
w^{k+1}=W(w^{k}+v^{k+1}-v^{k}),\qquad w^{0}=Wv^{0},\label{eq:gcon-2}
\end{equation}
which can be shown to also achieve consensus tracking. 

\begin{figure}
\subfloat[Centralized algorithm.\label{fig:calg}]{\begin{centering}
\includegraphics[width=1\columnwidth]{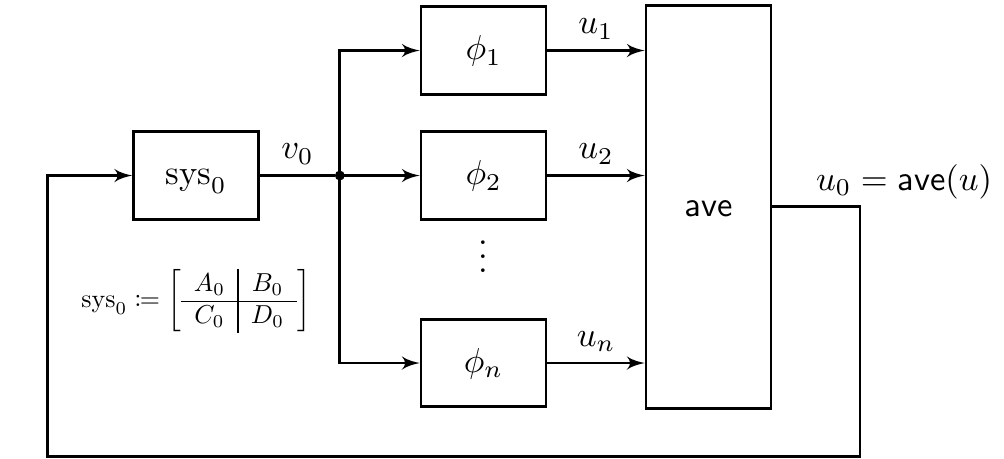}
\par\end{centering}
}

\subfloat[Corresponding decentralized algorithm.\label{fig:calg-de}]{\begin{centering}
\includegraphics[width=1\columnwidth]{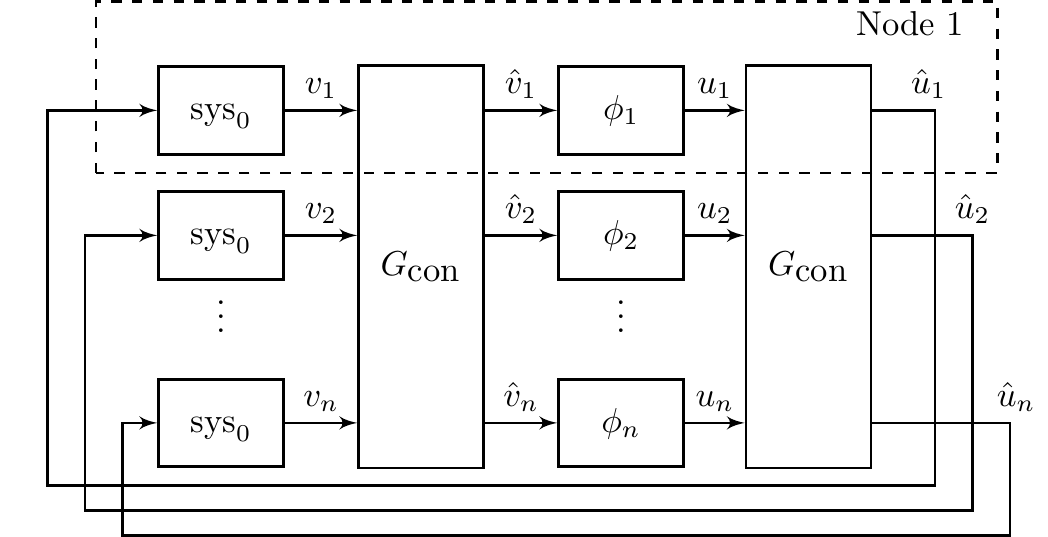}
\par\end{centering}
}

\caption{Comparison between a centralized algorithm and the corresponding decentralized algorithm given in Theorem\ \ref{thm:same_fp}.\label{fig:calg-compare}}
\end{figure}

The main idea behind converting an algorithm of the form\ (\ref{eq:calg}) or\ (\ref{eq:dalg}) into a decentralized one is to replace the $\ave$ operator with a system\ $\Gcon$ that achieves consensus tracking. In addition, the computation in\ (\ref{eq:calg-a}) and\ (\ref{eq:calg-b}) also needs to be decentralized, which can be handled by consensus tracking as well. For centralized algorithms of the form\ (\ref{eq:calg}), the resulting decentralized algorithm after conversion is described in the theorem below and also illustrated in Fig.\ \ref{fig:calg-compare}.
\begin{thm}
\label{thm:same_fp}Suppose $\Gcon$ is a system that achieves consensus tracking, and $\xi_{0}\opt$ is a possible steady-state value of $\xi_{0}$ in\ (\ref{eq:calg-a}). Then, $\xi_{0}\opt$ is also a steady-state value of $\xi_{i}$ in \begin{subequations}\label{eq:calg-gcon}
\begin{align}
\xi_{i}^{k+1} & =A_{0}\xi_{i}^{k}+B_{0}\hat{u}_{i}^{k}\label{eq:calg-gcon-a}\\
v_{i}^{k} & =C_{0}\xi_{i}^{k}+D_{0}\hat{u}_{i}^{k}\label{eq:calg-gcon-b}\\
u_{i}^{k} & =\phi_{i}(\hat{v}_{i}^{k})\label{eq:calg-gcon-c}
\end{align}
\end{subequations}for $i\in\{1,2,\dots,n\}$, where $\hat{u}=\Gcon u$ and $\hat{v}=\Gcon v$.
\end{thm}
\begin{IEEEproof}
Recall that if $u\opt$ is a possible steady-state value of $u$, then the corresponding steady-state value of $\hat{u}=\Gcon u$ is given by $\hat{u}\opt=Ju\opt$ or equivalently $\hat{u}_{i}\opt=\ave(u\opt)$. Suppose $(\xi_{0}\opt,v_{0}\opt,u\opt)$ is a steady-state value of $(\xi_{0},v_{0},u)$ in\ (\ref{eq:calg}). From\ (\ref{eq:calg-a}), we know $\xi_{0}\opt=A_{0}\xi_{0}\opt+B_{0}\ave(u\opt)$, which implies that\ (\ref{eq:calg-gcon-a}) is satisfied when $\xi_{i}^{k}=\xi_{0}\opt$ and $\hat{u}_{i}^{k}=\ave(u\opt)$ for all $k$. By checking\ (\ref{eq:calg-gcon-b}) and\ (\ref{eq:calg-gcon-c}) in a similar way, one can verify that $(\xi_{0}\opt,v_{0}\opt,u_{i}\opt,v_{0}\opt,\ave(u\opt))$ is a steady-value of $(\xi_{i},v_{i},u_{i},\hat{v}_{i},\hat{u}_{i})$ in\ (\ref{eq:calg-gcon}) for all $i\in\{1,2,\dots,n\}$. 
\end{IEEEproof}
\begin{rem}
\label{rem:same_fp_d}Theorem\ \ref{thm:same_fp} holds similarly for distributed algorithms of the form\ (\ref{eq:dalg}), in which case the corresponding decentralized algorithm becomes
\begin{align*}
\xi_{i}^{k+1} & =A_{i}\xi_{i}^{k}+B_{\text{loc},i}u_{i}^{k}+B_{i}\hat{u}_{i}^{k}\\
v_{i}^{k} & =C_{i}\xi_{i}^{k}+D_{\text{loc},i}u_{i}^{k}+D_{i}\hat{u}_{i}^{k}\\
u_{i}^{k} & =\phi_{i}(v_{i}^{k}),
\end{align*}
where $\hat{u}=\Gcon u$. We no longer need to apply consensus tracking on $v$ as in Theorem\ \ref{thm:same_fp}, because $\ave(u)$ in\ (\ref{eq:dalg-a}) and\ (\ref{eq:dalg-b}) is the only computation that prevents decentralization.
\end{rem}

Theorem\ \ref{thm:same_fp} and Remark\ \ref{rem:same_fp_d} give an equivalent decentralized algorithm, which admits the same steady state solution as the original algorithm. We give a few examples to illustrate how to use the result to decentralize an existing algorithm. We first give an example on the application to the gradient descent algorithm in Example\ \ref{exa:grad-ori}.
\begin{example}
[Decentralized gradient descent]\label{exa:grad-de} We use the system given in\ (\ref{eq:gcon-1}) as $\Gcon$, whose state space model can be written as
\[
\zeta^{k+1}=W\zeta^{k}+(W-I)v^{k},\qquad w^{k}=\zeta^{k}+v^{k},\qquad\zeta^{0}=0.
\]
Then, we can apply Theorem\ \ref{thm:same_fp} and obtain the following decentralized gradient descent algorithm:\begin{subequations}\label{eq:grad-de}
\begin{align}
\xi^{k+1} & =\xi^{k}-\eta\hat{u}^{k},\qquad v^{k}=\xi^{k},\qquad u^{k}=\nabla f(\hat{v}^{k})\label{eq:grad-de-a}\\
\zeta_{v}^{k+1} & =W\zeta_{v}^{k}+(W-I)v^{k},\qquad\hat{v}^{k}=\zeta_{v}^{k}+v^{k}\label{eq:grad-de-b}\\
\zeta_{u}^{k+1} & =W\zeta_{u}^{k}+(W-I)u^{k},\qquad\hat{u}^{k}=\zeta_{u}^{k}+u^{k},\label{eq:grad-de-c}
\end{align}
\end{subequations}where the initial condition is given by $\zeta_{v}^{0}=0$ and $\zeta_{u}^{0}=0$. Equation\ (\ref{eq:grad-de-a}) retains the original dynamics of the (centralized) gradient descent algorithm. The new equations\ (\ref{eq:grad-de-b}) and\ (\ref{eq:grad-de-c}) are due to the consensus tracking $\hat{v}=\Gcon v$ and $\hat{u}=\Gcon u$ of $v$ and $u$, respectively. It is not difficult to verify from\ (\ref{eq:grad-de})
\[
\hat{v}^{k+2}=2W\hat{v}^{k+1}-W^{2}\hat{v}^{k}-\eta(u^{k+1}-u^{k}),
\]
which recovers the dynamics of the DIGing algorithm in\ \cite{nedic_achieving_2017}. This should not come as a surprise, because the DIGing algorithm is based on consensus tracking of the average gradient. 
\end{example}
Next, we consider the ADMM algorithm in Example\ \ref{exa:admm-ori} and apply Remark\ \ref{rem:same_fp_d} to obtain a decentralized algorithm.
\begin{example}
[Decentralized ADMM]\label{exa:admm-de} We use the system given in\ (\ref{eq:gcon-2}) as $\Gcon$ and apply Remark\ \ref{rem:same_fp_d} to the ADMM algorithm in Example\ \ref{exa:admm-ori}. To make the algorithm more readable, we use the original optimization variables $x,y\in\reals^{n\times d}$ instead of $\xi$ and write the algorithm as
\begin{align*}
x^{k+1} & =Jx^{k}-(y^{k}+w^{k})/\rho,\qquad y^{k+1}=Jw^{k}-w^{k}\\
v^{k} & =Jx^{k}-(y^{k}+w^{k})/\rho,\qquad w^{k}=\nabla f(v^{k}).
\end{align*}
To apply Remark\ \ref{rem:same_fp_d}, we only need to replace $Jx$ with the consensus tracking $\hat{x}=\Gcon x$ of $x$ (and similarly for $Ju$). The corresponding decentralized algorithm is given by\begin{subequations}\label{eq:admm-de}
\begin{align}
x^{k+1} & =\hat{x}^{k}-(y^{k}+w^{k})/\rho,\qquad y^{k+1}=\hat{w}^{k}-w^{k}\label{eq:admm-de-a}\\
v^{k} & =\hat{x}^{k}-(y^{k}+w^{k})/\rho,\qquad w^{k}=\nabla f(v^{k})\label{eq:admm-de-b}\\
\zeta_{x}^{k+1} & =W\zeta_{x}^{k}+(W^{2}-W)x^{k},\qquad\hat{x}^{k}=\zeta_{x}^{k}+x^{k}\label{eq:admm-de-c}\\
\zeta_{w}^{k+1} & =W\zeta_{w}^{k}+(W^{2}-W)w^{k},\qquad\hat{w}^{k}=\zeta_{w}^{k}+w^{k},\label{eq:admm-de-d}
\end{align}
\end{subequations}where the initial condition is given by $\zeta_{x}^{0}=0$ and $\zeta_{w}^{0}=0$. Similar to Example\ \ref{exa:grad-de}, equations\ (\ref{eq:admm-de-a}) and\ (\ref{eq:admm-de-b}) retain the original ADMM dynamics, except that $Jx^{k}$ and $Jw^{k}$ are replaced respectively by $\hat{x}^{k}$ and $\hat{w}^{k}$, which are the output from consensus tracking given in\ (\ref{eq:admm-de-c}) and\ (\ref{eq:admm-de-d}).
\end{example}

\subsection{Discussions}

The result in Theorem\ \ref{thm:same_fp} can be viewed as a \emph{separation principle} for designing decentralized optimization algorithms. Specifically, a decentralized optimization algorithm can be formed by a non-decentralized base optimization algorithm (e.g., gradient descent, ADMM) and a decentralized consensus tracking system $\Gcon$. The system\ $\Gcon$ can be viewed as an approximation of the averaging operator $\ave$ that appears in the base algorithm; the faster $\Gcon$ reaches consensus, the better the approximation. The separation principle, however, does not require an explicit separation in time scale between the base algorithm and consensus tracking, which has been used in some previous work on decentralized gradient descent\ \cite{chen_fast_2012,jakovetic_fast_2014}. 

For converting centralized algorithms of the form\ (\ref{eq:calg}), an additional consensus tracking system\ $\Gcon$ is required. The conversion procedure can be interpreted as follows. Recall that the computation in\ (\ref{eq:calg-a})\textendash (\ref{eq:calg-b}) still takes place centrally (within ``$\mathrm{sys}_{0}$'' in Fig.\ \ref{fig:calg}), even when $\ave$ is replaced by a decentralized implementation. To make the computation decentralized, we create $n$ identical copies of $\mathrm{sys}_{0}$, one at each node. Despite being identical to each other, each $\mathrm{sys}_{0}$ will generate a different output $v_{i}$, because we can no longer guarantee $\hat{u}_{1}=\hat{u}_{2}=\cdots=\hat{u}_{n}$ after $\ave$ is replaced by $\Gcon$. Therefore, we need to use an additional $\Gcon$ to enforce consensus among all copies of $\mathrm{sys}_{0}$. 

To illustrate the importance of the additional $\Gcon$ for centralized algorithms, consider again the gradient descent algorithm in Example\ \ref{exa:grad-ori}. Without the additional $\Gcon$, the resulting decentralized algorithm would become\begin{subequations}\label{eq:grad-de-1}
\begin{align}
\xi^{k+1} & =\xi^{k}-\eta\hat{u}^{k},\qquad v^{k}=\xi^{k},\qquad u^{k}=\nabla f(v^{k})\label{eq:grad-de-a-1}\\
\zeta_{u}^{k+1} & =W\zeta_{u}^{k}+(W-I)u^{k},\qquad\hat{u}^{k}=\zeta_{u}^{k}+u^{k},\label{eq:grad-de-c-1}
\end{align}
\end{subequations}where the initial condition is given by $\zeta_{u}^{0}=0$. In steady state, we must have $\hat{u}\opt=0$ from\ (\ref{eq:grad-de-a-1}) and hence $\zeta_{u}\opt+u\opt=0$ from\ (\ref{eq:grad-de-c-1}). We also know $\vOne^{T}\zeta_{u}^{k}=0$ for all $k$ based on the initial condition and\ (\ref{eq:grad-de-c-1}). From these, we can only conclude $\vOne^{T}u\opt=0$ or equivalently $\sum_{i=1}^{n}f_{i}(v_{i}\opt)=0$. Therefore, we cannot obtain an optimal solution unless we have $v_{1}\opt=v_{2}\opt=\cdots=v_{n}\opt$, which could have been enforced by the additional $\Gcon$.

One benefit brought by the separation principle is that it allows us to derive a different decentralized algorithm by simply changing the consensus tracking system\ $\Gcon$. For example, using the system in\ (\ref{eq:gcon-2}) instead as $\Gcon$, we can obtain another decentralized gradient descent algorithm given by 
\begin{align*}
\xi^{k+1} & =\xi^{k}-\eta\hat{u}^{k},\qquad\hat{v}^{k}=\hat{\xi}^{k},\qquad u^{k}=\nabla f(\hat{v}^{k})\\
\zeta_{\xi}^{k+1} & =W\zeta_{\xi}^{k}+(W^{2}-W)\xi^{k},\qquad\hat{\xi}^{k}=\zeta_{\xi}^{k}+\xi^{k}\\
\zeta_{u}^{k+1} & =W\zeta_{u}^{k}+(W^{2}-W)u^{k},\qquad\hat{u}^{k}=\zeta_{u}^{k}+u^{k}.
\end{align*}
Another benefit of separation is reflected in the analysis of the resulting decentralized algorithm using the IQC framework proposed in\ \cite{lessard_analysis_2016}, which provides automated convergence analysis of optimization algorithms. Separation allows us to immediately reuse existing results derived for the base optimization algorithms, whereas we only need to incorporate the consensus tracking system\ $\Gcon$ into the IQC framework. This will be discussed in detail in Section\ \ref{sec:analysis}.

\section{Convergence Analysis\label{sec:analysis}}

We now show how to apply the IQC framework for automated convergence analysis of the decentralized algorithms obtained through Theorem\ \ref{thm:same_fp}. Throughout this section, we assume each $f_{i}$ is $\mu$\emph{-strongly convex} and \emph{$\beta$-smooth}, i.e., there exist $\mu>0$ and $\beta>0$ such that 
\[
\mu\norm{x-y}^{2}\leq\left(\nabla f_{i}(x)-\nabla f_{i}(y)\right)^{T}(x-y)\leq\beta\norm{x-y}^{2}.
\]
holds for all $x,y\in\mathbb{R}^{d}$. This assumption enables us to simplify the presentation and is not a limitation of the IQC analysis framework. For example, a similar IQC-based analysis has been developed when the assumption on strong convexity is removed\ \cite{hu_dissipativity_2017,fazlyab_analysis_2018,han_computational_2019}.

\subsection{IQC preliminaries\label{subsec:IQC-preliminaries}}

Many optimization algorithms, including the gradient descent method and ADMM presented in Examples\ \ref{exa:grad-ori} and\ \ref{exa:admm-ori}, can be viewed as a feedback interconnection of the form\begin{subequations}\label{eq:fbk_inter}
\begin{align}
\xi^{k+1} & =A\xi^{k}+Bu^{k},\qquad v^{k}=C\xi^{k}+Du^{k}\label{eq:lti}\\
u^{k} & =\phi(v^{k}).\label{eq:nl_map}
\end{align}
\end{subequations} We assume that the feedback connection is well-posed, which holds for both the gradient descent method and ADMM. Convergence analysis of an optimization algorithm becomes stability analysis of the interconnection\ (\ref{eq:fbk_inter}), which can be handled under the IQC framework. In the IQC framework, nonlinearity is treated as an uncertain system whose input $v$ and output $u$ are constrained by a quadratic inequality of the form
\begin{equation}
(z^{k}-z\opt)^{T}M(z^{k}-z\opt)\geq0,\quad\forall k,\label{eq:iqc}
\end{equation}
where $z=\Psi(v,u)$ is a ``filtered'' version of $v$ and $u$ given by
\[
\psi^{k+1}=A_{\Psi}\psi^{k}+B_{\Psi}\left[\begin{array}{c}
v^{k}\\
u^{k}
\end{array}\right],\quad z^{k}=C_{\Psi}\psi^{k}+D_{\Psi}\left[\begin{array}{c}
v^{k}\\
u^{k}
\end{array}\right].
\]
(Equation\ (\ref{eq:iqc}) is a special case of IQC called \emph{pointwise} IQC. Refer to\ \cite{lessard_analysis_2016} for more general IQCs.) For example, when $u^{k}=\phi(v^{k})=\nabla f(v^{k})$, based on strong convexity and smoothness of $f$, an IQC for $\phi$ is given by
\[
(z_{i}^{k}-\zopt_{i})^{T}\left(\left[\begin{array}{cc}
2\mu\beta & (\mu+\beta)\\
* & -2
\end{array}\right]\otimes I_{d}\right)(z_{i}^{k}-\zopt_{i})\geq0,
\]
where $z^{k}=(v^{k},u^{k})$. It has been shown in\ \cite{lessard_analysis_2016} that\ (\ref{eq:fbk_inter}) converges exponentially (linearly in the terminology used in optimization) with rate $\tau$ if there exists $P\succ0$ such that
\begin{align*}
 & \left[\begin{array}{cc}
A^{T}PA-\tau^{2}P & A^{T}PB\\
* & B^{T}PB
\end{array}\right]\otimes I_{d}\\
 & \quad+\left(*\right)^{T}(M\otimes I_{d})\left(\left[\begin{array}{cc}
C & D\\
0 & I
\end{array}\right]\otimes I_{d}\right)\preceq0,
\end{align*}
which is equivalent to 
\[
\left[\begin{array}{cc}
A^{T}PA-\tau^{2}P & A^{T}PB\\
* & B^{T}PB
\end{array}\right]+\left[*\right]^{T}M\left[\begin{array}{cc}
C & D\\
0 & I
\end{array}\right]\preceq0.
\]
The last step is called lossless dimensionality reduction\ \cite[Sec. 4.2]{lessard_analysis_2016} and has a useful interpretation: for convergence analysis, we can assume $d=1$ without loss of generality. This is consistent with the well-known fact that the convergence rate of many optimization algorithms does not depend on the dimension\ $d$ of the optimization variable. From here on, we will assume $d=1$ to simplify notation.

\subsection{Convergence analysis with a known gossip matrix\label{subsec:known_W}}

\begin{figure}[tbph]
\hfill{}\subfloat[Exact $W$ vs. unknown $W$ ($\sigma_{2}=0.2$, $W\in\reals^{2\times2}$).\label{fig:rate_admm_compare}]{\begin{centering}
\includegraphics[width=2.5in]{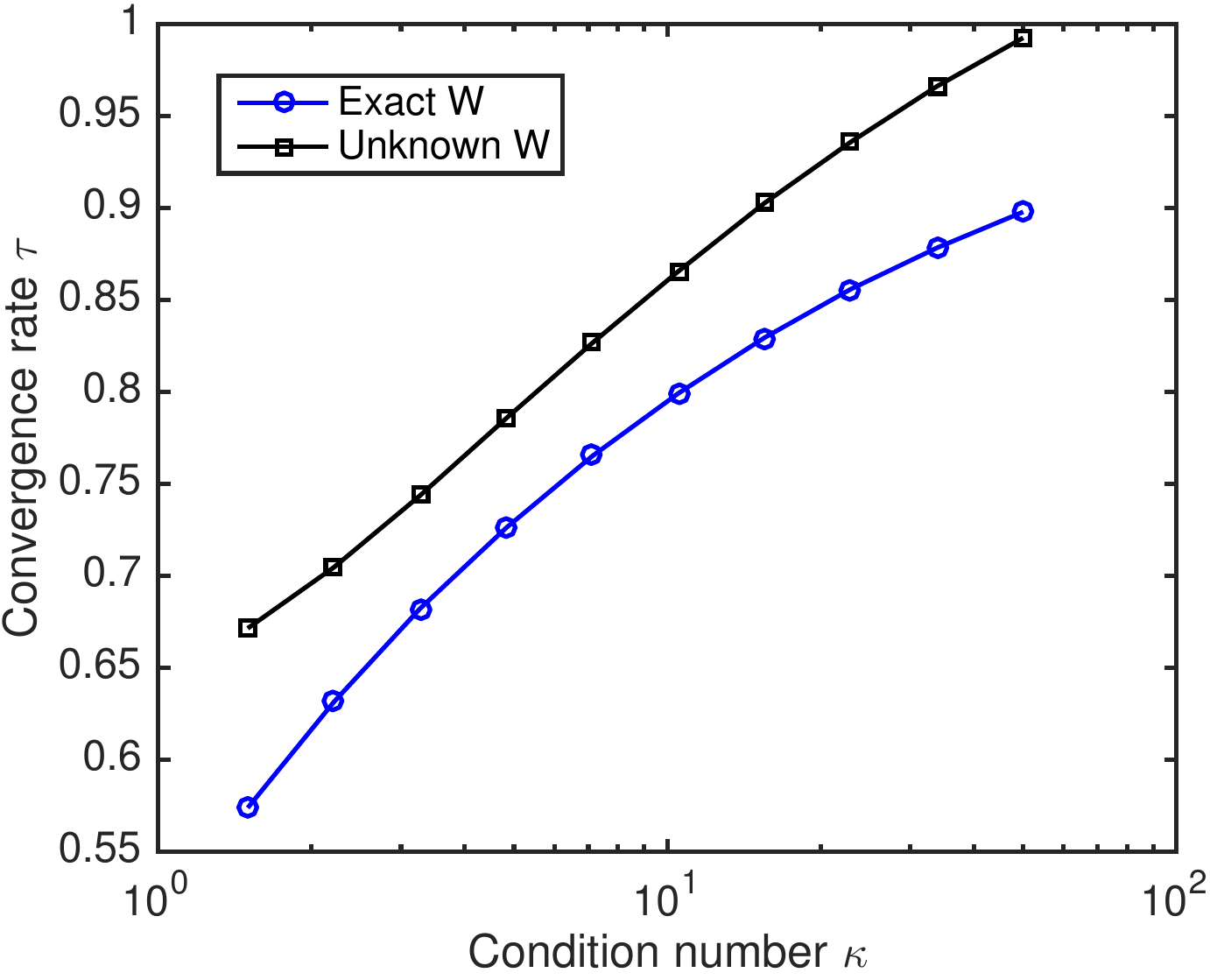}
\par\end{centering}
}\hfill{}

\hfill{}\subfloat[Worst-case convergence rate when only $\sigma_{2}$ is known.\label{fig:rate_admm_worst}]{\begin{centering}
\includegraphics[width=2.5in]{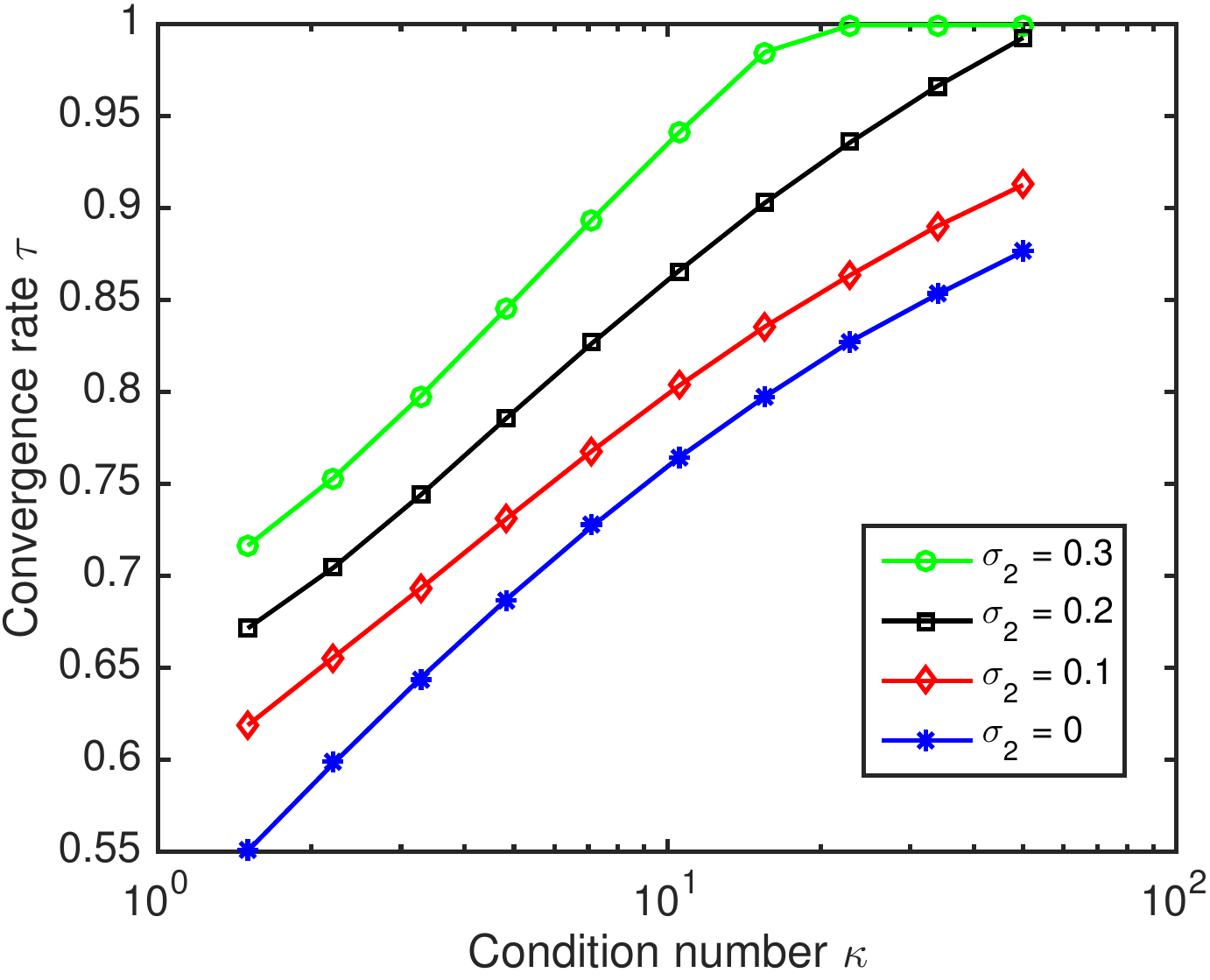}
\par\end{centering}
}\hfill{}

\caption{Convergence rate of decentralized ADMM ($\rho_{0}=1$).}

\end{figure}

We now use the IQC framework to verify the convergence of a decentralized algorithm obtained from Theorem\ \ref{thm:same_fp}. If the gossip matrix $W$ in $\Gcon$ is known, we can include the dynamics of\ $\Gcon$ in\ (\ref{eq:lti}) while keeping the same IQC for the nonlinear map in\ (\ref{eq:nl_map}). For example, we can write the decentralized ADMM algorithm\ (\ref{eq:admm-de}) in the form given in\ (\ref{eq:fbk_inter}) by choosing $\xi=(x,y,\zeta_{x},\zeta_{w})$. We adopt the normalization in\ \cite{nishihara_general_2015} and chose $\rho=\rho_{0}\sqrt{\mu\beta}$ for some fixed $\rho_{0}>0$ so that the convergence rate only depends on the condition number $\kappa\deq\beta/\mu$. The resulting system is given by
\[
\left[\begin{array}{c|c}
A & B\\
\hline C & D
\end{array}\right]=\left[\begin{array}{cccc|c}
W & -I & I & 0 & -I\\
0 & 0 & 0 & I & W-I\\
W^{2}-W & 0 & W & 0 & 0\\
0 & 0 & 0 & W & W^{2}-W\\
\hline W & -I & I & 0 & -I
\end{array}\right].
\]
The IQCs are given by $(z^{j,k}-z^{j,\star})^{T}M_{j}(z^{j,k}-z^{j,\star})\geq0$ for $j\in\{1,2,3\}$, where $z^{1,k}=(v^{k},w^{k})$, 
\begin{equation}
M_{1}=\left[\begin{array}{cc}
2\rho_{0}^{-2}I_{n} & \rho_{0}^{-1}(\kappa^{\frac{1}{2}}+\kappa^{-\frac{1}{2}})I_{n}\\
* & -2I_{n}
\end{array}\right],\label{eq:Mgrad_normalized}
\end{equation}
$z^{2,k}=\vOne^{T}\zeta_{x}^{k}$, $M_{2}=-1$, $z^{3,k}=\vOne^{T}\zeta_{w}^{k}$, and $M_{3}=-1.$ The IQC for $z^{1,k}$ comes from the properties of $f$. The other two IQCs encode the constraints $\vOne^{T}\zeta_{x}^{k}=0$ and $\vOne^{T}\zeta_{w}^{k}=0$, which are a result of the dynamics given by\ (\ref{eq:admm-de-c}) and\ (\ref{eq:admm-de-d}) under the zero initial condition $\zeta_{x}^{0}=0$ and $\zeta_{w}^{0}=0$. We computed the convergence rate $\tau$ for
\[
W=\left[\begin{array}{cc}
0.6 & 0.4\\
0.4 & 0.6
\end{array}\right].
\]
The result is shown in blue in Fig.\ \ref{fig:rate_admm_compare}. As expected, the convergence rate $\tau$ becomes slower as $\kappa$ increases. 

\subsection{Convergence analysis with an unknown gossip matrix}

The above IQC analysis relies on knowing the exact gossip matrix $W$. If only the second-largest singular value $\sigma_{2}\deq\sigma_{\max}(J-W)$ of $W$ is known, we are no longer able to include $\Gcon$ directly into the system dynamics\ (\ref{eq:lti}). Instead, we choose to treat $\Gcon$ as an uncertain system that can be characterized also by IQC. 

\begin{figure}[tbph]
\begin{centering}
\includegraphics[width=2.75in]{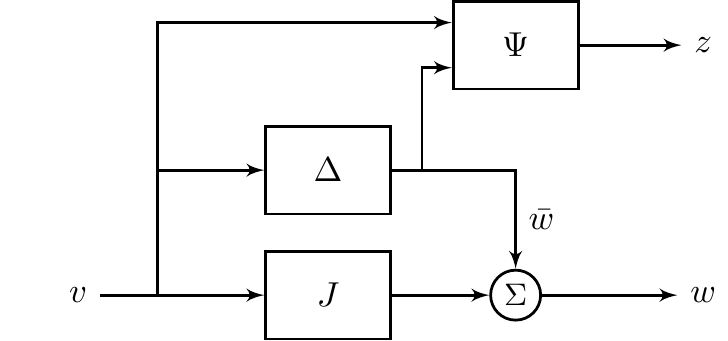}
\par\end{centering}
\caption{IQC characterization of $\Gcon$ with an unknown gossip matrix.\label{fig:gcon_iqc}}
\end{figure}

For the purpose of illustration, we will derive an IQC characterization (Fig.\ \ref{fig:gcon_iqc}) of the consensus tracking system\ $\Gcon$ given in\ (\ref{eq:gcon-2}), whose input is $v$ and output is $w$. Instead of constructing the filter $\Psi$ for $(v,w)$, we will construct the filter based on a different input-output pair in order to better capture certain important properties of\ (\ref{eq:gcon-2}). Because the steady-state value of $v$ and $w$ are related by $\wopt=J\vopt$, we define $\bar{w}\deq w-Jv$ so as to eliminate the steady-state component. Then, it can be shown that $\bar{w}$ and $v$ satisfy
\[
\bar{w}^{k+1}=W\bar{w}^{k}+(W-J)(v^{k+1}-v^{k}),\qquad\bar{w}^{0}=(W-J)v^{0}.
\]
As a result, we have $\vOne^{T}\bar{w}^{0}=0$ and $\vOne^{T}\bar{w}^{k+1}=\vOne^{T}\bar{w}^{k}$ for all $k$, i.e., $\vOne^{T}\bar{w}^{k}=0$ for all $k$. We write $W=J+\Jo\bar{W}\Jo$, where $\bar{W}\in\reals^{n\times n}$ satisfies $\sigma_{\max}(\bar{W})=\sigma_{2}$. Then, we have 
\begin{equation}
\bar{w}^{k+1}=\Jo\bar{W}\Jo(\bar{w}^{k}+v^{k+1}-v^{k}).\label{eq:gcon-2c}
\end{equation}
We define an uncertain system whose input $v$ and output $\bar{w}$ satisfy\ (\ref{eq:gcon-2c}). Notice 
\begin{equation}
\norm{\bar{w}^{k+1}}\leq\sigma_{2}\norm{\Jo(\bar{w}^{k}+v^{k+1}-v^{k})},\label{eq:sigma2_constr}
\end{equation}
which can be described by an IQC with $z^{1,k}\deq(\bar{w}^{k},\Jo(\bar{w}^{k}+v^{k}-v^{k-1}))$ and
\begin{equation}
M_{1}=\left[\begin{array}{cc}
-I & 0\\
* & \sigma_{2}^{2}I
\end{array}\right].\label{eq:M1_Psi}
\end{equation}
We also need to encode the constraint $\vOne^{T}\bar{w}^{k}=0$, which was not captured by\ (\ref{eq:sigma2_constr}). This can be done using an IQC with $z^{2,k}=\bar{w}^{k}$ and $M_{2}=-J$. The filter $\Psi$ that generates $z$ from $(v,\bar{w})$ is given by
\begin{equation}
\left[\begin{array}{c|c}
A_{\Psi} & B_{\Psi}\\
\hline C_{\Psi1} & D_{\Psi1}\\
\hdashline[2pt/2pt]C_{\Psi2} & D_{\Psi2}
\end{array}\right]=\left[\begin{array}{cc|cc}
0 & 0 & I & 0\\
0 & 0 & 0 & I\\
\hline 0 & 0 & 0 & I\\
-\Jo & \Jo & \Jo & 0\\
\hdashline[2pt/2pt]0 & 0 & 0 & I
\end{array}\right].\label{eq:Psi_def}
\end{equation}

We can carry out the convergence analysis by forming an interconnection of the following four systems: $(v,x)=\Gadmm(w,\hat{x},\hat{w})$, where $\Gadmm$ is given by\ (\ref{eq:admm-de-a})\textendash (\ref{eq:admm-de-b}) as
\[
\left[\begin{array}{cc|ccc}
0 & -I & I & I & 0\\
0 & 0 & -I & 0 & I\\
\hline 0 & -I & I & I & 0\\
I & 0 & 0 & 0 & 0
\end{array}\right],
\]
$z_{\nabla}=(v,w)$, $z_{x}=\Psi(x,\hat{x}-Jx)$, and $z_{w}=\Psi(w,\hat{w}-Jw$), where $\Psi$ is given in\ (\ref{eq:Psi_def}). The input to the interconnection is $(w,\hat{x},\hat{w})$, and the output from the interconnection is $(z_{\nabla},z_{x},z_{w})$. The output obeys the following IQCs: $z_{\nabla}$ is constrained by $M_{1}$ in\ (\ref{eq:Mgrad_normalized}), and both $z_{x}$ and $z_{w}$ are constrained by the two matrices $M_{1}$ and $M_{2}$ associated with $\Psi$. 

We computed the worst-case convergence rate when only $\sigma_{2}$ is known. We applied another dimensionality reduction introduced in\ \cite{sundararajan_robust_2017} so that the convergence rate is independent of\ $n$. The result for different values of $\sigma_{2}$ is shown in Fig.\ \ref{fig:rate_admm_worst}. As expected, the convergence rate becomes slower as $\sigma_{2}$ increases. As can be seen from Fig.\ \ref{fig:rate_admm_compare}, the worst-case convergence analysis is more conservative than the result from Sec.\ \ref{subsec:known_W}. (When $W$ is $2\times2$, it is uniquely determined by $\sigma_{2}$.) However, we would like to emphasize that our purpose is to demonstrate the capability of integrating consensus tracking into the IQC framework. Moreover, it is also possible to reduce the level of conservatism by enriching the class of IQCs (e.g., Zames-Falb IQC) used in analysis.

\section{Conclusions}

We have proposed a separation principle for designing decentralized algorithms used in consensus optimization. Specifically, a decentralized optimization algorithm can be constructed by combining a non-decentralized base optimization algorithm and decentralized consensus tracking; the latter replaces the averaging operation that appears in the base algorithm. The separation principle provides modularity in both the design and analysis of algorithms. For design, the principle allows one to choose any combination of base algorithm and consensus tracking algorithm. For analysis, modularity is enabled by the automated convergence analysis based on IQC, which is capable of integrating consensus tracking, regardless of whether the underlying gossip matrix is known. As a result, convergence of the decentralized algorithm can be readily verified as long as the base algorithm already has an existing IQC characterization; the computation is as simple as calculating the interconnection of multiple linear dynamical systems coming from the base algorithm and consensus tracking. The workflow of design and analysis has been illustrated using a decentralized ADMM algorithm. We believe that the same principle also applies to other optimization problems that only require local information sharing, e.g., when locally coupled objective function and/or constraints are present (cf.\ \cite[Sec. 7.2]{boyd_distributed_2011}).

\bibliographystyle{abbrv}
\bibliography{CDC2019}

\end{document}